\documentclass[11pt,reqno]{amsart}

\usepackage{amsfonts,amsthm, amscd, epsfig, amsmath, amssymb,enumerate,psfrag}

\numberwithin{equation}{section}

\DeclareMathOperator{\Spec}{Spec}   

\DeclareMathSymbol{\leqslant}{\mathalpha}{AMSa}{"36} % nicer `smaller or equal' 
\DeclareMathSymbol{\geqslant}{\mathalpha}{AMSa}{"3E} % nicer `larger or equal' 
\DeclareMathSymbol{\eset}{\mathalpha}{AMSb}{"3F}     % nicer `emptyset' 
\renewcommand{\leq}{\;\leqslant\;}                   % redef. of < or = 
\renewcommand{\geq}{\;\geqslant\;}                   % redef. of > or = 
\newcommand{\grad}{\nabla} 

\newcommand{\la}{\label} 
 
\newcommand{\be}{\begin{equation}}

\def\1{\ifmmode {1\hskip -3pt \rm{I}} \else {\hbox {$1\hskip -3pt \rm{I}$}}\fi}

\newtheorem{Th}{Theorem}[section] 
\newtheorem{Le}[Th]{Lemma} 
\newtheorem{Pro}[Th]{Proposition}

\newcommand{\cE}{\ensuremath{\mathcal E}}

\newcommand{\cH}{\ensuremath{\mathcal H}}

\newcommand{\cL}{\ensuremath{\mathcal L}}

\newcommand{\cX}{\ensuremath{\mathcal X}}

\newcommand{\bbR}{{\ensuremath{\mathbb R}} }

\newcommand{\bbZ}{{\ensuremath{\mathbb Z}} }

\newcommand{\var}{{\rm Var} }

\newcommand{\wt}{\widetilde}

\let\a=\alpha      \let\e=\varepsilon
 \let\g=\gamma       \let\l=\lambda
            
  \let\s=\sigma    
  \let\z=\zeta

\def\\{\hfill\break}
\def\thsp{\thinspace}

\def\?{\mskip -10mu}
\def\tc{\thsp | \thsp}

\newcommand{\LRW}{\cL^{RW}}
\newcommand{\lRW}{\l^{RW}}
\newcommand{\SRW}{S^{RW}}
\newcommand{\LIC}{\cL^{I\!P}}
\newcommand{\lIC}{\l^{I\!P}}
\newcommand{\mIC}{\mu^{I\!P}}
\newcommand{\SIC}{S^{I\!P}}
\newcommand{\LEX}{\cL^{EP}}
\newcommand{\lEX}{\l^{EP}}
\newcommand{\SEX}{S^{EP}}
\newcommand{\LCEP}{\cL^{CEP}}
\newcommand{\lCEP}{\l^{CEP}}
\newcommand{\SCEP}{S^{CEP}}
\newcommand{\LCY}{\cL^{CP}}
\newcommand{\lCY}{\l^{CP}}
\newcommand{\SCY}{S^{CP}}
\newcommand{\LMP}{\cL^{M\!P}}
\newcommand{\lMP}{\l^{M\!P}}

\newcommand{\SMP}{S^{M\!P}}

\usepackage{graphicx}

\begin{document}

\title[Proof of Aldous' spectral gap conjecture]
{Proof of Aldous' spectral gap conjecture}

\author[P. Caputo]{Pietro Caputo}
\address{%
	Pietro Caputo \hfill\break
	\indent Dipartimento di Matematica , Universita' di Roma Tre, Italy {\em and} \hfill\break
	\indent Department of Mathematics, University of California, Los Angeles, USA}
\email{caputo\@@mat.uniroma3.it}

\author[T.M. Liggett]{Thomas M.\ Liggett} 
\address{%
	Thomas M.\ Liggett \hfill\break
	\indent  Department of Mathematics, University of California, Los Angeles, USA} 
\email{tml\@@math.ucla.edu}

\author[T. Richthammer]{Thomas Richthammer} \address{Thomas Richthammer \hfill\break
 	\indent Department of Mathematics, University of California, Los Angeles, USA}
\email{richthammer\@@math.ucla.edu}

%    \subjclass is required.


%\date{Last revised: September 26, 2009}

%    Abstract is required.

\begin{abstract}
Aldous' spectral gap conjecture asserts that on any 
graph the random walk process and the 
random transposition (or interchange) process have the same spectral gap. 
We prove the conjecture using a recursive 
strategy. 
The approach is a natural extension of the  
method already used to prove 
the validity of the conjecture on trees.
The novelty is an idea based on electric 
network reduction, which reduces 
the problem to the proof of an explicit inequality 
for a random transposition operator involving both positive and 
negative rates. 
The proof of the latter inequality uses suitable coset decompositions 
of the associated matrices with rows and columns indexed by permutations. 
\end{abstract}

\maketitle

\thispagestyle{empty}

\section{Introduction}

\noindent
Spectral gap analysis 
plays an important role in the study  
of the convergence to equilibrium of reversible Markov chains. 
We begin by reviewing some well known facts about Markov chains and 
their spectra. For more details we refer to \cite{AldousFill}.

\subsection{Finite state, continuous time Markov chains.} 
\label{secmc}
Let us consider a continuous time Markov chain $Z = (Z_t)_{t \geq 0}$ with 
finite state space $S$ and transition rates $(q_{i,j}: i \neq j \in S)$
such that $q_{i,j} \geq 0$. We will always assume that the Markov chain is irreducible and satisfies 
$$
q_{i,j} = q_{j,i} \quad \text{ for all }i \neq j.
$$ 
Such a Markov chain is reversible with respect to the uniform distribution $\nu$ on $S$, which is the unique stationary distribution of the chain. 
The infinitesimal generator $\cL$ of the Markov chain is
defined by 
$$
\cL g (i) = \sum_{j\in S} q_{i,j} (g(j)-g(i))\,,
$$
where $g:S\to\bbR$ and $i\in S$. 
The matrix corresponding to the linear operator $\cL$ is 
the transition matrix $Q = (q_{i,j})_{i,j}$, 
where $q_{i,i} := - \sum_{j \neq i} q_{i,j}$, 
and the corresponding quadratic form is 
$$
\sum_{i \in S} g(i) \cL g(i) = \sum_{i,j\in S} q_{i,j} g(i)(g(j)-g(i)) 
= - \frac12\sum_{i,j\in S} q_{i,j} (g(j)-g(i))^2\,. 
$$
Thus, $-\cL$ is positive semi-definite and symmetric, which implies 
that its spectrum is of the form 
$\Spec (-\cL) = \{ \l_i: 0 \leq i \leq |S|-1\}$, where  
$$
0 = \l_0 < \l_1 \leq \ldots \leq \l_{|S|-1}.
$$
The spectral gap $\l_1$ is characterized as the largest constant $\l$ such that 
\begin{equation}\la{gap1}
\frac12\sum_{i,j\in S} q_{i,j} (g(j)-g(i))^2 \geq \l \sum_{i \in S}g(i)^2\,
\end{equation}
for all $g:V\to\bbR$ with $\sum_{i}g(i)=0$. The significance of 
$\l_1$ is its interpretation as the asymptotic rate of convergence to the 
stationary distribution: 
$$
P_i(Z_t = j) = \nu(\{j\}) + a_{i,j}e^{-\l_1 t} +o(e^{-\l_1 t}) \qquad 
\text{ for } t \to \infty, 
$$
where typically $a_{i,j} \neq 0$ (and more precisely  
$a_{i,i} >0$ for some $i$). For this reason $\frac{1}{\l_1}$ is 
sometimes referred to as the relaxation time of the Markov chain, 
and it is desirable to have an effective way of calculating $\l_1$. 
Aldous' conjecture relates the spectral gap of the random walk 
on a finite graph to that of more complicated processes on the same graph. 
This can be very important in applications, 
since generally speaking it is easier to compute or estimate 
(e.g.\ via isoperimetric inequalities)
the spectral gap of the random walk than that of the other processes considered, which have much larger state spaces.\\ 

We say that the Markov chain with state space $S_2$ and generator $\cL_2$
is a sub-process of the chain with state space $S_1$ and generator $\cL_1$ 
if there is a contraction of $S_1$ onto $S_2$, i.e.\ 
if there is a surjective map $\pi: S_1 \to S_2$ such that 
\begin{equation}\la{mcprojection}
\cL_1(f \circ \pi) = (\cL_2 f) \circ \pi \quad 
\text{ for all $f:S_2 \to \bbR$.}
\end{equation}
In this case, suppose that $f$ is an eigenfunction of $-\cL_2$ 
with eigenvalue $\l$. Then  
$-\cL_1(f \circ \pi) = (-\cL_2 f) \circ \pi = \l f \circ \pi$ 
and $f \circ \pi \neq 0$ for $f \neq 0$, so 
$f\circ \pi$ is an eigenfunction of $-\cL_1$ with the same eigenvalue $\l$. 
Thus,
$$
\Spec(-\cL_2) \subset \Spec(-\cL_1), 
$$
and, in particular, the spectral gap of the first process is smaller than or equal 
to that of the second process. Identity (\ref{mcprojection}) is an
example of a so-called intertwining relation; see e.g.\
\cite{Diaconis-Fill} for more details on such relations and their 
applications.

\subsection{Random walk and interchange process on a weighted graph.} \la{secprocesses}
%In the following subsections we will define different stochastic processes 
%on a finite graph $G$, all of which are Markov chains of the type considered 
%in Section~\ref{secmc}. 
Let $G = (V,E)$ be an undirected complete graph on $n$ vertices; without loss of generality we assume that its vertex set 
is $V =\{1,\dots,n\}$. Furthermore $G$ is a weighted graph 
in that we are given a collection of edge weights (or conductances) 
$c_{xy} \geq 0$, for $xy = \{x,y\} \in E$. Since we want 
the processes defined below to be irreducible,  
we will assume that the skeleton graph, 
i.e.\ the set of edges $xy$ where $c_{xy}>0$, is connected. 
If we want to stress the dependence of one of the processes described
below on the underlying weighted graph, we will write $\cL(G)$ and 
$\l_1(G)$ for its generator and gap.  
Finally, we note that considering complete graphs only 
is no loss of generality, 
since edges with weight 0 can be thought of as being ``absent''.

\subsubsection{Random walk}
\begin{figure}[htb]
\psfrag{1}{$1$}
\psfrag{2}{$2$}
\psfrag{3}{$3$}
\psfrag{4}{$4$}
\psfrag{5}{$5$}
\includegraphics[scale = .55]{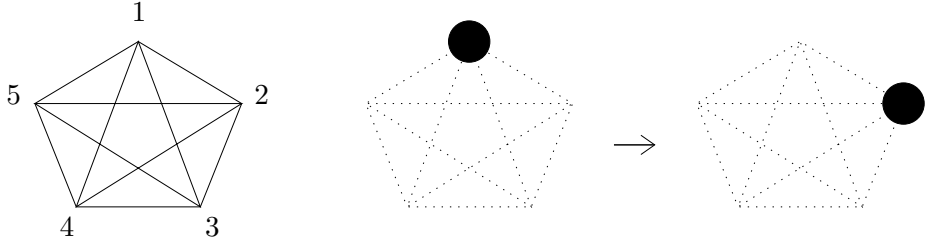}
\caption{Random walk on $V=\{1,2,3,4,5\}$. 
The picture shows the underlying graph and a transition
from state 1 to state 2.}  
\label{figrw}
\end{figure}

The (1-particle) random walk on $G$ is the Markov chain
in which a single particle jumps from vertex 
$x \in  V$ to $y \neq x$ at rate $c_{xy}$; see Figure~\ref{figrw}. 
Formally, its state space is $\SRW = V = \{1,2,\ldots,n\}$ and 
its generator is defined by 
$$
\LRW f (x) = \sum_{y \neq x} c_{xy}(f(y)-f(x))\, , \quad \text{ for } 
f: V \to \bbR, \; x \in V.
$$
By Section~\ref{secmc}, $-\LRW$ has $n=|\SRW|$ nonnegative eigenvalues 
and a positive spectral gap $\lRW_1 >0$.

\subsubsection{Interchange process} 
\begin{figure}[htb]
\psfrag{1}{$1$}
\psfrag{2}{$2$}
\psfrag{3}{$3$}
\psfrag{4}{$4$}
\psfrag{5}{$5$}
\psfrag{A}{$1$}
\psfrag{B}{$2$}
\psfrag{C}{$3$}
\psfrag{D}{$4$}
\psfrag{E}{$5$}
\includegraphics[scale = .55]{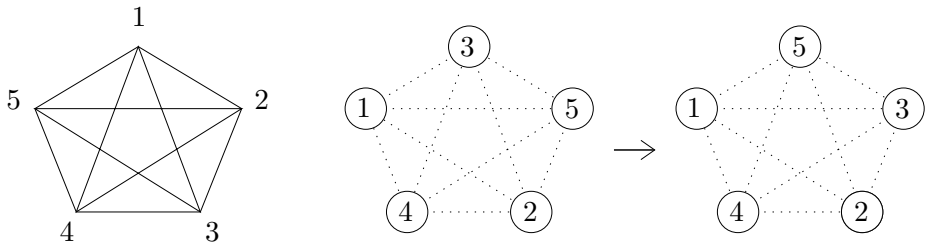}
\caption[]{Interchange process on $V=\{1,2,3,4,5\}$. 
The picture shows the underlying graph and a transition from 
state $\eta = \left(\begin{smallmatrix}
1 & 2 & 3 & 4 & 5\\ 5 & 3 & 1 & 4 & 2 \end{smallmatrix} \right)$ 
to
$\eta^{1,2} = \eta \tau_{1,2} = \left( \protect\begin{smallmatrix}
1 & 2 & 3 & 4 & 5\\ 5 & 3 & 2 & 4 & 1 \end{smallmatrix} \right)$. 
Note that in this notation the first row refers to the labels. 
}
\label{figinterchange}
\end{figure}

\la{secinterchange}
In the interchange process, a state is an assignment of $n$ labeled particles 
to the vertices of $G$ in such a way that each vertex is occupied by exactly one particle. The transition from a state $\eta$ to a state $\eta^{xy}$ 
(occurring with rate $c_{xy}$) interchanges the particles at vertices 
$x$ and $y$; see Figure~\ref{figinterchange}.
For a formal definition, let $\cX_n$ denote the set of permutations of $V=\{1,\dots,n\}$, 
and for $\eta \in \cX_n$ and $xy \in E$ let 
$\eta^{xy} = \eta \tau_{xy}$, where $\tau_{xy} \in \cX_n$ is
the transposition of $x$ and $y$. The interchange process on $G$ 
is the Markov chain with state space $\SIC = \cX_n$ and 
generator
$$
\LIC f (\eta) = \sum_{xy \in E} c_{xy} (f(\eta^{xy}) - f(\eta)), \quad 
\text{ where } f: \SIC \to \bbR, \,\eta \in \SIC.
$$
We use $\eta_x$ to denote the label of the particle at 
$x$, while $\xi_i = \xi_i(\eta)$ will be used to denote the position of the particle labeled $i$. By Section~\ref{secmc}, $-\LIC$ has $|\SIC| = n!$ nonnegative eigenvalues and a positive spectral gap $\lIC_1 >0$.  
The random walk can be obtained as a sub-process of the interchange process by 
ignoring all particles apart from the one with label 1; more precisely 
the map $\pi: \SIC \to \SRW$, $\pi(\eta) := \xi_1(\eta)$ is a contraction
in the sense of \eqref{mcprojection}. 
Thus, 
$$
\Spec(-\LRW) \subset \Spec(-\LIC), 
$$
and in particular, 
\be\la{ineq}
\lIC_1 \leq \lRW_1 \,.
\end{equation} 
%Aldous' conjecture 
\subsection{Main result}
Our main result states that inequality (\ref{ineq}) is, in fact, an equality: 
\begin{Th}\la{main}
For all weighted graphs $G$, the interchange process and the random walk 
have the same spectral gap:
\begin{equation}\la{th}
\lIC_1(G)=\lRW_1(G)\,.
\end{equation}
\end{Th}
A weaker form of Theorem \ref{main} involving only unweighted graphs 
had been conjectured by Aldous around 1992, and since then it has been mostly
referred to as Aldous' spectral gap conjecture in the literature. 
Related observations can be found in Diaconis and Shahshahani's paper \cite{DS}, and in 
the comparison theory developed
by Diaconis and Saloff-Coste \cite{Diaconis-Saloff}.
%
%\begin{Con}[Aldous 1992]\la{aldous}
%For all weighted graphs $G$, the interchange process and the random walk 
%have the same spectral gap:
%$$\lIC_1(G)=\lRW_1(G)\,.$$ 
%\end{Con}
%
%The purpose of this paper is to prove this conjecture. 
%Before discussing the conjecture and outlining our strategy for its proof
%we will describe its consequences for other processes on $G$. 

The problem has received a lot of attention in recent years 
- the conjecture was stated as an open problem on David Aldous' web page 
\cite{Aldous} and in the influential 
monographs \cite{AldousFill,LPW}. In the    
meantime, various special cases have been obtained.  
The first class of graphs that was shown to satisfy the conjecture 
is the class of unweighted complete 
graphs (i.e.\ $c_{xy} = 1$ for all $xy \in E$). 
Diaconis and Shahshahani computed all eigenvalues of the interchange process 
in this case using the irreducible representations 
of the symmetric group \cite{DS}. 
Similar results were obtained for unweighted star graphs in \cite{FOW}. 
Recently, remarkable work of Cesi pushed this
algebraic approach further to obtain the conjecture 
for all unweighted complete multipartite graphs 
\cite{F}. 

An alternative approach based on recursion 
was proposed by Handjani and Jungreis \cite{HJ} (see also Koma and Nachtergaele \cite{KN} for similar results) 
who proved the conjecture for all weighted trees. 
The same ideas were recently used by Conomos and Starr \cite{CS}, and Morris \cite{Morris}, 
to obtain an asymptotic version of the conjecture for boxes in the lattice $\bbZ^d$ with unweighted edges.  The basic recursive approach in \cite{HJ} 
has been recently rephrased in purely algebraic terms, 
see \cite[Lemma 3.1]{F2}.

In order to prove Theorem \ref{main} we develop a general recursive 
approach based on the idea of network reduction; see Section \ref{recu}. 
The method, inspired by the theory of resistive networks, 
allows us to reduce the proof of the theorem to the proof 
of an interesting comparison inequality for random  
transposition operators on different weighted graphs; 
see Theorem \ref{leoctopus}.
%see Theorem \ref{main} and Theorem \ref{leoctopus}. 

After a 
preliminary version \cite{CLR1} of this paper appeared, we learned that 
the same recursive strategy had been discovered around the same time
independently, and from a slightly different perspective,
by Dieker \cite{Dieker}. The comparison inequality alluded to above 
%claim in our Theorem \ref{leoctopus} 
was conjectured to hold in 
both \cite{CLR1} and \cite{Dieker}. 

The comparison inequality will be proved in Section \ref{secoctopus}.
The main idea for this proof is a decomposition of 
the associated matrix into a covariance matrix 
and a correction matrix (a Schur complement). 
A delicate analysis based on block 
decompositions corresponding to suitable cosets of the 
permutation group reveals that the 
correction matrix is nonnegative definite.

Some immediate consequences of Theorem \ref{main} for other natural Markov
chains 
associated to finite weighted graphs are discussed in Section \ref{sec_rem}.

We end this introductory section with a collection of 
known properties of the spectrum of 
the interchange process that can be deduced from the algebraic approach. 
We refer to 
\cite{DS,FOW,F} and references therein for more details. 
These facts are not needed 
in what follows and the 
reader may safely jump to the next section. However, 
we feel that 
the algebraic point of view provides a natural decomposition of the 
spectrum that is worth mentioning. 

\subsection{Structure of the spectrum of $-\LIC$}
In Section~\ref{secinterchange} we saw that 
$\Spec(-\LRW) \subset \Spec(-\LIC)$.
One can go a little further and show that, if $0=\lRW_0<\lRW_1\leq\cdots\leq \lRW_{n-1}$ 
are the eigenvalues of $-\LRW$, then for 
$k\geq 0$ and $1\leq i_1<\cdots<i_k\leq n-1$, 
\begin{equation}\la{sumo}
\lRW_{i_1}+\cdots+\lRW_{i_k} \in \Spec(-\LIC)\,.
\end{equation}
The corresponding eigenfunction is the  
antisymmetric product of the $k$ one-particle eigenfunctions of $\lRW_{i_1},\dots,\lRW_{i_k}$.
In particular, the eigenvalue 
\begin{equation}\la{sums}
\lRW_1 + \cdots + \lRW_{n-1}
={\rm Tr}(-\LRW) = 2 \sum_{xy \in E}c_{xy} \,,
\end{equation}
is associated with 
functions that are antisymmetric in all particles, i.e.\ 
multiples of the alternating function 
$h(\eta)={\rm sign}(\eta)$. 
(This also follows directly from $h(\eta^{xy})-h(\eta) = -2h(\eta)$.)
From the representation theory of the symmetric group  one can compute 
(see below) the 
multiplicity of all eigenvalues of the form (\ref{sumo}),
and one finds that the overwhelming majority 
(for large $n$) of the spectrum of $-\LIC$ are not of this form.

The vector space of functions $f:\cX_n\to\bbR$ is equivalent to a direct sum 
$\oplus_{\a}\cH_\a$, where $\a$ ranges over all (equivalence classes of the) irreducible representations of the symmetric group.
Since the latter are in one to one correspondence with the partitions 
of $n$, one can identify $\a$ with a Young diagram $\a=(\a_1,\a_2,\dots)$, where the $\a_i$
form a non-increasing sequence of nonnegative integers such that $\sum_{i}\a_i=n$. 

Each subspace $\cH_\a$ is in turn a 
direct sum $\cH_\a=\oplus_{j=1}^{d_\a}\cH^{j}_\a$,
of subspaces $\cH_\a^j$, each of dimension $d_\a$, where the positive integer 
$d_\a$ is the dimension of 
the irreducible representation $\a$. In particular, the 
numbers $d_\a$ satisfy $\sum_{\a} (d_\a)^2=n!$.
The subspaces $\cH_\a^j$ are invariant for 
the action of the generator $-\LIC$, so that $-\LIC$ can be diagonalized within each $\cH^j_\a$. 
Subspace $\cH^i_\a$ will produce $d_\a$ eigenvalues 
$\l_k(\a)$, $k=1,\dots,d_\a$. 
Some of these may coincide if the weights have suitable symmetries 
(for instance, if $G$ is the complete graph with $c_{xy} = 1$
for all $xy \in E$, 
then they all coincide and 
$-\LIC$ is a multiple of the identity matrix in each $\cH_\a$, cf.\ \cite{DS}) but 
in the general weighted case they will be distinct. 
On the other hand, for a given $\a$, the 
eigenvalues coming from $\cH^i_\a$ are identical to those coming from $\cH^j_\a$, for all $i,j=1,\dots,d_\a$, 
so that each eigenvalue $\l_k(\a)$ 
will appear with multiplicity $d_\a$ in the spectrum of $-\LIC$. Moreover, using known 
expressions for the characters of transpositions, one can compute explicitly
the sum $$\sum_{k=1}^{d_\a}\l_k(\a)\,,$$ for every irreducible representation $\a$, as a 
function of the edge weights.  For instance, when $\a$ is the partition
$(n-1,1,0,\dots)$, which has $d_\a = n-1$, one obtains the relation (\ref{sums}).
The trivial partition $(n,0,\dots)$ has dimension $1$ and the only eigenvalue is $0$. This is the space of constant functions. Similarly, the alternating partition 
$(1^n,0\dots)$ ($n$ ones and then all zeros), 
has dimension $1$ and the only eigenvalue is $2 \sum_{xy \in E} c_{xy}$. 
It can be shown that the eigenvalues of the form \eqref{sumo} come from the L-shaped partitions $\a=(n-k,1^k,0,\dots)$, each with dimension $d_\a=\binom{n-1}k$. So the total number of eigenvalues of the form 
\eqref{sumo} is $\sum_{k=0}^{n-1} \binom{n-1}k ^2 = \binom{2(n-1)}{n-1}$.

Finally, using known relationships between conjugate irreducible representations, 
see e.g.\ \cite[2.1.8]{JamesKerber}, \cite[(2.12)]{F2},
one can show that  
the spectrum of $-\LIC$ can be decomposed into pairs of 
eigenvalues $\l,\l'$ such that 
$$\l+\l'=2 \sum_{xy \in E}c_{xy}\,,$$ 
where $\l,\l'$ are associated with conjugate Young diagrams.

\section{A recursive approach based on network reduction}\la{recu}

\begin{figure}[h]
\centerline{
\psfrag{x}{{\Small $5$}}\psfrag{a}{{\Small $1$}}\psfrag{b}{{\Small $2$}}\psfrag{c}{{\Small $3$}}\psfrag{e}{{\Small $4$}}
\psfrag{arr}{$\Longrightarrow$}
\psfig{file=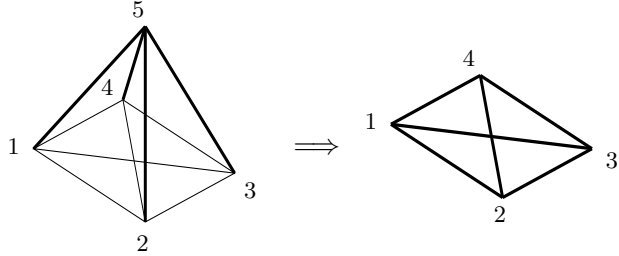,height=1.3in, width=3.2in}
}
\caption{Reduction 
of a 5-vertex graph to a 4-vertex graph at $x=5$.}  
\label{fig:red54}
\end{figure}

\noindent
Given a weighted graph $G=(V,E)$ as above and a point $x\in V$ we consider the 
{\em reduced network} obtained by removing the vertex $x$. 
This gives a new graph $G_x$ with vertex set $V_x:= V\setminus\{x\}$, 
edge set $E_x = \{yz \in E: y,z \neq x\}$
and edge conductances $\wt c_{yz}\geq c_{yz}$ defined by 
\begin{equation}\la{redrates}
\wt c_{yz} = c_{yz} + c^{*,x}_{yz}\,,\quad\; c^{*,x}_{yz}:=\frac{c_{xy}c_{xz}}
{\sum_{w\in V_x} c_{xw}}\,,
\end{equation}
for $yz\in E_x$. We refer to $G_x$ as the reduction of $G$ at $x$ 
or simply as the reduced graph at $x$. 
This is the general version of more familiar network reductions such as 
series resistance (from 3 to 2 vertices) or star--triangle transformations (from 4 to 3 vertices); see Figure~\ref{fig:red54}. 
We refer to 
\cite{DoyleSnell,LyonsPeres,AldousFill} for the classical 
probabilistic point of view on electric networks. 

\medskip

\subsection{Random walk on the reduced network}
We first show that the spectral gap of the random walk 
on the reduced network is not smaller than the original random 
walk spectral gap:
\begin{Pro}\la{redlemma}
The spectral gaps of the random walks on a weighted graph $G$ and 
the corresponding reduced graph $G_x$ satisfy  
$$\lRW_1(G_x)\geq \lRW_1(G)\,.$$
\end{Pro}
\proof
We will use the shorthand notation $L = \LRW(G)$ and $L_x = \LRW(G_x)$ 
for the generators of the two random walks.
We first note that, for any graph $G$, 
$\lRW_1(G)$ can be characterized as the largest
constant $\l$ such that
\begin{equation}\la{gap2}
\sum_{z\in V} (L g(z))^2
\geq -\l\sum_{z\in V} g(z) L g(z)\,
\end{equation} 
holds for all $g:V\to \bbR$. 
To see this, observe that, for any $g,h: V\to\bbR$, 
$\sum_{z\in V}h(z)Lg(z) = \sum_{z\in V}g(z)Lh(z)$.
Thus, taking $h=Lg$, the left hand side of (\ref{gap2}) 
coincides with the quadratic form 
$$
\sum_{z\in V} g(z) L^2g(z)\,,
$$
and (\ref{gap2}) says that $L^2 + \l  L$ is nonnegative definite. 
Taking a basis which makes $L$ diagonal one sees that this holds iff $\l\leq \lRW_1(G)$. 

To prove the proposition, 
take a function $g:V\to\bbR$ harmonic at $x$, i.e.\ such that $Lg(x) = 0$. Then 
\begin{equation}\la{harm}
g(x) = \frac{\sum_{y\in V_x} c_{xy} g(y)}{\sum_{w\in V_x} c_{xw}} \,.
\end{equation}
For any $z\in V_x$, from (\ref{harm}) we have 
\begin{align*}
Lg (z) &= \sum_{y\in V_x} c_{zy}[g(y)-g(z)] + c_{zx}[g(x)-g(z)]
\\
& = \sum_{y\in V_x}  \left(c_{zy} + c_{zy}^{*,x}\right)[g(y)-g(z)]\,.
\end{align*}
In other words
$$
Lg (z) = \begin{cases} L_xg (z) & z\in V_x\\
0 & z=x.
\end{cases}
$$
%Therefore $\sum_{z\in V} g(z) Lg(z) = \sum_{z\in V_x}  g(z) L_x g(z)$ and 
%$\sum_{z\in V} Lg(z) Lg(z) = \sum_{z\in V_x}  L_xg(z) L_x g(z)$.
Applying (\ref{gap2}) to this function
we have 
\begin{align*}
\sum_{z\in V_x}  (L_xg(z))^2 &=\sum_{z\in V}  (Lg(z))^2\geq -\lRW_1(G) 
\sum_{z\in V}  g(z)Lg(z) \\
&= -\lRW_1(G)\sum_{z\in V_x}  g(z) L_x g(z)\,.
\end{align*}
Since the function $g$ is arbitrary on $V_x$, 
using \eqref{gap2} again, this time for the graph $G_x$,
we obtain $\l_1^{RW}(G_x)\geq \l_1^{RW}(G)$.
\qed

\bigskip

Proposition \ref{redlemma} 
generalizes the observation in \cite{HJ}
that if $G$ is a graph with a vertex $x$ of degree $1$ 
(i.e.\ only one edge out of $x$ has positive weight),
then the spectral gap of the random walk cannot decrease 
when we cancel $x$ and remove the only 
edge connecting it to the rest of $G$.
(In that case $\wt c_{yz} = c_{yz}$ since $x$ has degree $1$.) 

We end this subsection with a side remark on further 
relations between the generators $L = \LRW(G)$ and $L_x = \LRW(G_x)$.
When we remove a vertex, it is interesting to compare the 
energy corresponding to the removed branches with the energy coming from the 
new conductances. The following identity can be obtained with a straightforward  
computation. 
\begin{Le}
\la{harmo}
For any fixed $x\in V$ and any $g:V\to \bbR$, 
$$%\begin{equation}\la{hbo}
\sum_{y \in V_x} c_{xy} [g(y)-g(x)]^2 = \sum_{yz \in E_x} c^{*,x}_{yz} [g(y)-g(z)]^2
+ \frac{1}{\sum \limits_{y\neq x} c_{xy}}\,(Lg(x))^2
\,. 
$$
%\end{equation}
\end{Le}
Consider the operator $\wt L_x$ defined by $\wt L_x g(x) = 0$ and 
$\wt L_x g(z) = L_x g(z)$ for $z\neq x$, where $g:V\to \bbR$. Then $\wt L_x$ is the generator of the random walk on $G_x\cup \{x\}$, where $x$ is an isolated vertex.
Lemma \ref{harmo} implies that %, in the vector space $\bbR^V$, 
the quadratic form of $-\wt L_x$ is dominated by the quadratic form of $-L$. It follows from the Courant-Fisher min-max theorem that, if $\wt \l_0\leq \cdots\leq \wt \l_{n-1}$
denote the eigenvalues of $-\wt L_x$, 
then $\wt\l_i \leq \l_i^{RW}(G)$, $i=0,\dots,n-1$. 
Note that this is not in contradiction with the
result in Proposition \ref{redlemma} since, due to the isolated 
vertex $x$, one has $\wt \l_0=\wt \l_1=0$, and $\wt \l_{k+1}=\l^{RW}_k(G_x)$, $k=1,\dots,n-2$.
%the spectral gap 
%$\l_1^{RW}(G_x)$ coincides with $\wt\l_2$. 
While the bound in Proposition \ref{redlemma} will be sufficient for our purposes, it is worth pointing out that, as observed in \cite{Dieker}, at this point
standard results on interlacings can be used to prove
the stronger statement
$$
\l_j^{RW}(G)\leq \l_j^{RW}(G_x) \leq \l_{j+1}^{RW}(G)
\,,\quad j=1,\dots,n-2\,.
$$

\medskip

\subsection{Octopus inequality}
The following theorem summarizes the main technical ingredient 
we shall need. Here $\nu$ is the uniform probability measure 
on all permutations $\cX_n$, and we use the notation $\nu[f] = \int f\,d\nu$. 
The gradient $\grad$ is 
defined by 
$$
\grad_{xy} f(\eta) = f(\eta^{xy}) - f(\eta)\,.
$$
\begin{Th}\la{leoctopus}
For any weighted graph $G$ on $|V|=n$ vertices, 
for every $x \in V$ 
and $f:\cX_n\to\bbR$: 
\begin{equation}\la{octopus}
\sum_{y \in V_x} c_{xy}\, \nu[(\grad_{xy} f)^2] 
\,\geq  \sum_{yz \in E_x} c^{*,x}_{yz} \,\nu[(\grad_{yz} f)^2] \,.
\end{equation}
\end{Th}

\medskip

\noindent
Note that if $f(\eta)=g(\xi_1)$ is a function of one particle, then a simple computation gives
$$ \nu[(\grad_{uv} f)^2] = \frac2{n}(g(u)-g(v))^2\,,\quad\;uv\in E \,,$$ 
so that this special case of Theorem \ref{leoctopus} is contained in Lemma \ref{harmo}.
The identity in Lemma \ref{harmo} also shows that in this case the inequality is saturated by 
functions that are harmonic at $x$. 
On the other hand, the general case represents a nontrivial comparison inequality
between a weighted star graph and its complement, with weights
defined by (\ref{redrates}). 
Inspired by its tentacular nature 
we refer to the bound (\ref{octopus}) as the {\em octopus} inequality. 
We will give a proof 
of Theorem \ref{leoctopus} in Section~\ref{secoctopus}.

\subsection{Reformulation of the conjecture}
We shall use the following convenient notation: 
As above let $\nu$ denote the uniform probability measure $\cX_n$, 
$\grad$ the gradient and $b$ a generic edge, whose weight is denoted $c_b$. 
In this way $\LIC=\sum_b c_b \grad_b$ and the Dirichlet form 
$-\nu[f\LIC f]$ is 
$$
\cE(f) = \frac12\sum_b c_b \,\nu[(\grad_bf)^2]\,.
$$
The spectral gap $\lIC_1$ is the best constant $\l$ so that 
for all $f:\cX_n\to\bbR$:
\begin{equation}\la{gaps1}
\cE(f)\geq \l\,\var_{\nu}(f)\,,
\end{equation}
where $\var_{\nu}(f) = \nu[f^2]-\nu[f]^2$ is the variance of $f$ w.r.t.\ $\nu$.
In order to get some hold on the eigenvalues of the interchange process 
that are not eigenvalues of the random walk we introduce the vector space
\[
\begin{split}
\cH &= \{f:\cX_n\to\bbR : \nu[f\tc \xi_i]=0 \text{ for all } i \in V\}\\
&= \{f:\cX_n\to\bbR : \nu[f\tc \eta_x]=0 \text{ for all } x \in V\}, 
\end{split}
\]
where $\nu[\cdot\tc \xi_i]$ and $\nu[\cdot\tc \eta_x]$ are 
the conditional expectations given the position of the particle labeled $i$
and given the label of the particle at $x$ respectively. 
The equality in the definition of $\cH$ is a consequence of 
$$
\nu[\cdot\tc \xi_i](\eta)=\nu[\cdot\tc \xi_i = x] = \nu[\cdot\tc \eta_x = i]
= \nu[\cdot\tc \eta_x](\eta) \,,
$$
where $\eta\in\cX_n$ is  such that $\xi_i(\eta)=x$. Note that for
every $i$:
$$
\nu[\LIC f \tc \xi_i=x] = \sum_{y\neq x} c_{xy} \left( \nu[f \tc \xi_i=y]
- \nu[f \tc \xi_i = x ]\right)\,,
%\LRW(\nu[f \tc \xi_i]) \quad 
%\text{ for all } f: \cX_n \to \bbR,  
$$
for all $f: \cX_n \to \bbR$, and $x\in V$. In particular, 
$\cH$ is an invariant subspace for $-\LIC$, and if $f \notin \cH$ 
is an eigenfunction of $-\LIC$ with eigenvalue $\l$, 
then $\nu [f\tc \xi_i]\neq 0$ for some $i$, and 
$\nu [f\tc \xi_i=x]$, $x\in V$, is an eigenfunction of $-\LRW$ with the same eigenvalue $\l$.
It follows that 
$\cH$ contains all eigenfunctions corresponding to eigenvalues in 
$\Spec(-\LIC) \setminus \Spec(-\LRW)$.
Therefore, if $\mIC_1(G)$ denotes the 
smallest eigenvalue of $-\LIC$ 
associated to functions in $\cH$
(i.e.\ the best constant $\l$ in (\ref{gaps1}) 
restricting to functions $f\in\cH$), then for every graph $G$ one has
$$\lIC_1(G)=\min\{\lRW_1(G),\mIC_1(G)\}\,.$$
The assertion $\lIC_1(G)=\lRW_1(G)$ of Theorem~\ref{main} becomes then equivalent to 
\begin{equation}\la{lambdamu}
\mIC_1(G)\geq \lRW_1(G) \,.  
\end{equation}
In the rest of this section we show how the network reduction idea, assuming the validity of Theorem \ref{leoctopus},
yields a proof of Theorem \ref{main}.

\medskip

\subsection{Proof of Theorem \ref{main}}

We use the notation from the previous sections. In particular 
we write $\lRW_1(G_x)$ and $\lIC_1(G_x)$ for the spectral gaps of 
the random walk and the interchange process in the network reduced at 
$x$.
Let us first show that Theorem \ref{leoctopus} implies an estimate 
of $\mIC_1(G)$. 

\begin{Pro} \la{promumax}
For an arbitrary weighted graph $G$
\begin{equation}\la{bombom}
\mIC_1(G)\geq \max_{x\in V} \lIC_1(G_x)\,.
\end{equation}
\end{Pro}

\proof Let $f\in\cH$ and $x \in V$. 
Since $\nu [f\tc \eta_x]=0$,
we have $$\nu[f^2]=\var_\nu(f) = \nu[\var_\nu(f\tc\eta_x)]\,,$$ 
where $ \var_\nu(f\tc\eta_x)$
is the variance w.r.t.\ $\nu[\cdot\tc\eta_x]$. 
For a fixed value of $\eta_x$, $\nu[\cdot\tc\eta_x]$
is the uniform measure on the permutations on $V_x=V\setminus \{x\}$. 
Therefore using the spectral gap bound \eqref{gaps1} on the graph $G_x$ we have
$$
\lIC_1(G_x)\,\var_\nu(f\tc\eta_x) 
\leq \,\frac12\sum_{b: \,b\not\ni x}(c_{b} + c_b^{*,x})
\nu[(\grad_b f)^2\tc \eta_x]\,,
$$ 
with $c_b^{*,x}$ defined by (\ref{redrates}).
%Taking an average with some probability $q$ on $V$ and 
Taking the $\nu$-expectation we obtain:
$$
\lIC_1(G_x)\,\nu[f^2] \leq 
\frac12
\sum_{b: \,b\not\ni x} (c_b + c_b^{*,x})
\nu[(\grad_b f)^2]\,.
$$
From Theorem \ref{leoctopus}:
$$
\sum_{b: \,b\not\ni x}c_b^{*,x}\, \nu[(\grad_b f)^2]
\leq \sum_{b: \,b\ni x} c_b \,\nu[(\grad_b f)^2]\,.
$$ 
Therefore, 
\begin{equation}\la{bom}
\lIC_1(G_x) \,\nu[f^2]\leq \cE(f)\,.
\end{equation}
Since $x\in V$ and $f\in\cH$ were arbitrary, 
this proves that, for every $x\in V$, $\mIC_1(G)\geq \lIC_1(G_x)$, establishing the inequality \eqref{bombom}. \qed

\medskip

Propositions \ref{redlemma} and \ref{promumax}
allow us to conclude the proof by recursion. 
Indeed, note that $\lIC_1(G)=\lRW_1(G)$ is trivially true when $G=b$ 
is a single weighted edge $b$. (When $n=2$, the 
random walk and the interchange process 
are the same 2-state Markov chain.) 
If $G$ is a weighted graph on $n$ vertices, we assume 
that $\lIC_1(G')=\lRW_1(G')$ holds on every weighted graph $G'$ 
with $n-1$ vertices, in particular on $G_x$. Therefore 
$$
\mIC_1(G) \geq \max_{x\in V} \lIC_1(G_x) =\max_{x\in V} \lRW_1(G_x)
\geq \lRW_1(G),
$$
where we also have used Propositions \ref{redlemma} and \ref{promumax}.
Thus we have shown \eqref{lambdamu}, which is equivalent to 
$\lIC_1(G) = \lRW_1(G)$.

\section{Proof of the octopus inequality}\label{secoctopus}

For the proof of Theorem~\ref{leoctopus} we slightly change our notation as follows. We set 
$V=\{0,1,\ldots,n-1\}$ and $x=0$. The only rates appearing in \eqref{octopus} 
are $c_{0i}$, so we set 
\begin{align*}
&c_i := c_{0i} \;\;\text{ for }\, 1 \leq i\,, \;\;\quad 
c_0 := -\sum_{i \geq 1}  c_i \quad \\
&\quad\;\text{ and } \quad 
c := \sum_{1 \leq i \leq n-1} c_i^2 + \sum_{1 \leq i<j \leq n-1} c_ic_j\,.
\end{align*}
%in order to simplify notation. 
Note that $c_0 < 0$ and %we have 
\begin{equation}\la{relcoeff}
\sum_{i \geq 0} c_i = 0, \quad 
c = -\sum_{0 \leq i<j} c_ic_j
\quad \text{ and } \quad 
c_{ij}^{*,0} = - \frac{c_ic_j}{c_0}.
\end{equation}
Using this shorthand notation the octopus inequality \eqref{octopus} simplifies to 
\begin{equation}\la{octopus2}
-\sum_{0 \leq i<j} c_ic_j \sum_{\eta} (f(\eta \tau_{ij}) - f(\eta))^2 \geq 0,
\end{equation}
where $\tau_{ij}$ denotes the transposition of $i,j \in V$, i.e.\ $\eta \tau_{ij}=\eta^{ij}$. 
Thus it suffices to show that the matrix $C$ defined by 
\begin{equation}\la{defC}
C_{\eta, \eta'} = \left\{
\begin{array}{ll} 
c & \text{ if } \eta  = \eta'\\
c_ic_j & \text{ if } \eta \tau_{ij} = \eta'\\
0 & \text{ otherwise,}
\end{array}
\right.
\end{equation}
is positive semi-definite for every $n$ and all rates $c_1,\ldots,c_{n-1} \geq 0$.

\subsection{Decomposition of the matrix $C$}\la{secred}
In the following we write $A \geq B$ if the same inequality holds for the 
corresponding quadratic forms, i.e.\ if $A-B$ is positive semi-definite. 
Obviously, this defines a partial order %on symmetric matrices 
and 
we will repeatedly use the following simple facts 
for %symmetric 
square matrices $A,B$ and a real number $a$: 
\begin{align*}
&A \geq 0\,, \;B \geq 0 \;\Rightarrow\; A+B \geq 0\,; \qquad 
\;A \geq 0\,, \;a \geq 0 \;\Rightarrow aA \geq 0\,; \\
&
\qquad \qquad \qquad
\left( \begin{array}{cc} A &  0 \\  0 & B \end{array} \right)  \geq 0 
\;\,\Leftrightarrow\;\, A,B \geq 0\,.
\end{align*}
Note that 
every transposition takes even to odd permutations and vice versa, 
so $C$ has the block structure 
$$
C = \left( \begin{array}{cc} cI & X^t \\ X & cI 
\end{array}
\right), \quad \text{ where $I$ is the identity 
matrix},
$$
and we have used a basis which lists first all even permutations, 
and then all odd permutations. We have 
$$ 
\tilde{C} := 
\left( \begin{array}{cc}  \frac 1 c X^t X  &  X^t \\  
X & cI \end{array} \right) 
= 
\left( \begin{array}{cc}  \frac 1 {\sqrt{c}} X & \sqrt{c}I \end{array} \right)^t
\left( \begin{array}{cc}  \frac 1 {\sqrt{c}} X & \sqrt{c}I 
\end{array}
\right) \geq 0, 
$$
since $A^tA \geq 0$ for any matrix $A$, and $C$ and $\tilde{C}$ only differ by 
$$
C-\tilde{C} = \left( \begin{array}{cc}  \frac 1 c C' & 0 \\ 0 & 0 
\end{array} \right), \quad 
\text{ where } 
C' = c^2I - X^t X.
$$
$C'$ is a symmetric 
$\frac{n!}{2} \times \frac{n!}{2}$-matrix, to be referred to as the correction matrix. 
It coincides with $c$ times 
the Schur complement of the odd-odd 
block of $C$. The matrices $C$, $C'$ and $X$ only depend on the rates $c_1,\ldots,c_{n-1}$ and the  system size $n = |V|$, 
and whenever we want to stress this dependence 
we will write $C(n)$, $C'(n)$ and $X(n)$. 
By the above, the proof of Theorem \ref{leoctopus} 
will be complete once we show that $C'$ is positive semi-definite: 
\begin{equation}\la{cor}
C'(n) \geq 0 \quad \text{ for all } n \geq 2.  
\end{equation}

\subsection{Structure of the correction matrix} \la{seccor} 
It turns out that the correction matrix has a relatively simple structure: It can be written as a linear combination of matrices 
where the coefficients are products of rates and the matrices 
do not depend on the rates  at all. 
\begin{Le} \la{lecorrdecomp}
We have $C'(2) = 0$, $C'(3) = 0$ and 
\begin{equation}\la{corrdecomp}
C'(n) = \sum_{J \subset V: |J| = 4} -c_J A^J(n) 
\qquad \text{ for all } n \geq 4, 
\end{equation}
where $c_J := \prod_{i \in J} c_i$ and 
$A^J(n)$ is defined by
\begin{equation}\la{defmat}
A^J_{\eta, \eta'}(n) = \left\{
\begin{array}{ll} 
2 & \text{ if } \eta  = \eta'\\
2 & \text{ if $\eta^{-1}\eta'$ is a product of 2 disjoint 2-cycles}\\
& \hspace*{3.5 cm} \text{  with entries from  $J$}\\
-1 & \text{ if $\eta^{-1}\eta'$ is a 3-cycle with entries from $J$}\\
0 & \text{ otherwise}
\end{array}
\right.
\end{equation}
for all even permutations $\eta,\eta' \in \cX_V$.
\end{Le}

\proof 
We simply calculate $C'_{\eta,\eta'}$ for all even permutations 
$\eta,\eta'$ using $C' = c^2I - X^tX$. For $n=2$, 
$c = c_1^2$ and $X(2)$ is the $1\times 1$-matrix $X(2) = (-c_1^2)$, so 
$C'(2) = 0$. For $n=3$, $c=c_1^2+c_1c_2+c_2^2$ and 
$$X(3) = \left( \begin{array}{ccc}
c_1c_2 & -c_2(c_1+c_2) & -c_1(c_1+c_2) \\
-c_1(c_1+c_2)&c_1c_2 & -c_2(c_1+c_2)\\
-c_2(c_1+c_2) & -c_1(c_1+c_2) & c_1c_2
\end{array}
\right),$$
where the rows are indexed by the odd permutations 
$(12),(01),(02)$ and the columns are indexed by the even permutations 
$id,(021),(012)$ in that order. This gives 
$C'(3) = c^2 I(3) - X^t(3) X(3) = 0$. 

For $n \geq 4$ we observe that  $X_{\eta_1,\eta_2}(n) = 0$ unless $\eta_1$ and $\eta_2$ differ by a single 
transposition. Thus $C'_{\eta,\eta'}(n) = 0$ unless $\eta$ and $\eta'$ 
differ by a product of exactly two transpositions. Note that such a product 
of two transpositions can be a product of 
two disjoint transpositions (i.e.\ 2-cycles), a 3-cycle,  
or the identity. 
\\
(a) If $\eta^{-1} \eta'$ is a product of two disjoint 2-cycles, 
e.g.\ $(0 1)(2 3)$, 
a complete list of decompositions of  $\eta^{-1} \eta'$ into a product 
of two transpositions is $(0 1)(2 3) = (2 3)(0 1)$, so using 
$K := \{0,1,2,3\}$ we have 
$$
C'_{\eta,\eta'}(n) = -(c_0c_1c_2c_3 + c_2c_3c_0c_1) 
= 2 (-c_K)\,.
$$
(b) If $\eta^{-1} \eta'$ is a 3-cycle, e.g.\ $(0 1 2)$, a complete list 
of decompositions of $\eta^{-1} \eta'$ into a product of two transpositions 
is $(0 1 2) = (01)(20) = (12)(01) = (20)(12)$, so using $K := \{0,1,2\}$ 
we have 
\begin{align*}
C'_{\eta,\eta'}(n) &= -(c_0c_1c_2c_0 + c_1c_2c_0c_1 + c_2c_0c_1c_2)
\\
&= c_K \sum_{i \notin K} c_i
= (-1) \sum_{J \supset K,|J|=4} -c_J\,.
\end{align*}
(c) If $\eta^{-1} \eta' = id$, we have $\eta^{-1} \eta' = \tau^2$ for 
every transposition $\tau$, so we have 
$$
C'_{\eta,\eta'}(n) = c^2 - \sum_{i<j}(c_ic_j)^2 = 2 \sum_{J:|J|=4} -c_J. 
$$
Here we have used 
\begin{align*}
&\Big(\sum_{i<j} c_ic_j\Big)^2 - \sum_{i<j} (c_ic_j)^2 \\
&\qquad= 2 \sum_{i<j<k} \Big(c_i^2c_jc_k  + c_ic_j^2c_k + c_ic_jc_k^2\Big)  
+  6 \sum_{i <j < k < l} c_ic_jc_kc_l\\
&\qquad=2 \Big(\sum_{i<j<k} c_ic_jc_k \sum_l c_l  -  4 \sum_{i <j<k < l} c_ic_jc_kc_l\Big)+ 
6 \sum_{i < j<k < l} c_ic_jc_kc_l \\
&\qquad= - 2 \sum_{i <j<k < l} c_ic_jc_kc_l.
\end{align*}
Thus we have checked \eqref{corrdecomp} entrywise. 
\qed

\medskip

We already know that $C'(2) = 0$ and $C'(3) = 0$. 
In order to motivate the following lemmata let us also look at 
$C'(4)$ and $C'(5)$:
Using the shorthand notation 
$$
A := A^{\{0,1,2,3\}}(4)\,, \quad \text{ and } \;\;  
A^{(i)} := A^{\{0,1,2,3,4\}\setminus\{i\}}(5)\,,\quad \text{ for } 0 \leq i \leq 4,
$$
the decomposition \eqref{corrdecomp} of the correction matrices gives 
\[
\begin{split}
C'(4) &= - c_0c_1c_2c_3 A\,, \quad\text{ and } \\
C'(5) &=  -c_0c_2c_3c_4A^{(1)}- \ldots 
-c_0c_1c_2c_3A^{(4)} - c_1c_2c_3c_4A^{(0)}.
\end{split}
\]
For $C'(4)$ we observe that $ -c_0c_1c_2c_3 \geq 0$, so it suffices 
to show that $A \geq 0$. For $C'(5)$ we observe that 
$$
-c_0 c_2c_3c_4 = (c_1+c_2+c_3+c_4) c_2c_3c_4 \geq c_1c_2c_3c_4,
$$
and similarly for the coefficients of $A^{(2)}$, $A^{(3)}$ and $A^{(4)}$. 
If $A^{(i)} \geq 0$, this implies 
$$
C'(5) \geq c_1c_2c_3c_4 ( A^{(1)} + A^{(2)} + A^{(3)} + A^{(4)} -
A^{(0)}),   
$$
and since $c_1c_2c_3c_4 \geq 0$ we are done once we have shown  
that the matrix in the parentheses is $\ge0$. 
For general $n$ we will need the following two lemmata. In their 
proofs we will repeatedly use the notation $\cX_K$ and $\cX_K^+$ 
for the set of all permutations on a set $K$ and the set 
of all even permutations on $K$. 

\begin{Le} \la{lepropA}
For all $n \geq 4$ and $J \subset V$ with $|J|=4$: 
\begin{equation}\la{propA}
A^J(n) \geq 0\,.
\end{equation} 
\end{Le}

\proof 
Consider the block structure of $A^J(n)$ corresponding to the blocks 
formed by the $n!/4!$ left cosets of $\cX_J^+$ in $\cX_V^+$. 
By definition of $A^J(n)$ in \eqref{defmat}, the diagonal block 
corresponding to the coset $\cX_J^+$ can be identified with 
$A := A^{\{0,1,2,3\}}(4)$ 
(if $J$ is identified with $\{0,1,2,3\})$. Furthermore
$A^J_{\eta,\eta'}(n)$ only depends on $\eta^{-1}\eta'$, and thus
$A^J_{\eta \s,\eta \s'}(n) = A^J_{\s,\s'}(n)$ for all $\s,\s' \in \cX_J^+$
and $\eta \in \cX_V^+$,
which implies that all diagonal blocks of $A^J(n)$ are equal, 
and thus they are copies of $A$.  
Finally, $A^J_{\eta,\eta'}(n) =0$ unless $\eta^{-1} \eta' \in \cX_J^+$, 
which shows that all non-diagonal blocks of $A^J(n)$ are $0$.
Because of this block decomposition of $A^J(n)$ we only have to show $A \geq 0$.

By \eqref{defmat} $A = A^{\{0,1,2,3\}}(4)$ is a symmetric 
$12 \times 12$ matrix with entries $2,-1,0$ only. 
Using a computer algebra program one can check that $A$ has the eigenvalues 0 (with multiplicity 10) and 12 
(with multiplicity 2), which implies the assertion for $n=4$. 
For the sake of completeness we will also show how to obtain the spectrum 
of $A$ without using a computer: 
The matrix $A$ is indexed by $\cX_4^+ := \cX^+_{\{0,1,2,3\}}$. 
We note that $\cX_4^+$ consists of 
the identity, 3 permutations that are a product of two 
disjoint 2-cycles and 8 permutations that are 3-cycles. 
Furthermore 
$H := \{id, (01)(23), (02)(13), (03)(12)\}$ is a subgroup of $\cX_4^+$,
and we consider the decomposition of $A$ into blocks corresponding
to the 3 left cosets of $H$ in $\cX_4^+$. Two permutations from the same 
coset $\eta H$ differ by an element of $H$, whereas two 
permutations from different cosets can't differ by an element of $H$, 
i.e.\ they have to differ by a 3-cycle. 
Thus by \eqref{defmat} $A$ has the block structure
\begin{equation} \la{A4}
A = \left( 
\begin{array}{ccc}
2 E_4 & (-1) E_4 & (-1) E_4 \\
(-1) E_4 & 2 E_4 & (-1) E_4\\ 
(-1) E_4 & (-1) E_4 & 2 E_4 
\end{array}\right)
= 3 \left( \begin{array}{ccc}
 E_4 & 0 & 0 \\
0 &  E_4 & 0\\ 
0 & 0 &  E_4 
\end{array}\right)
- E_{12},
\end{equation}
where $E_m \in \bbR^{m \times m}$ is the matrix with all entries equal to $1$. 
The matrix $E_n$ has the eigenvalues $0$ (with multiplicity $n-1$)
and $n$ (with multiplicity 1), and the eigenvector corresponding to the 
eigenvalue $n$ is $(1,\ldots,1)$. Furthermore the two matrices in the above
decomposition of $A$ commute. This implies that 
$A$ has the eigenvalues $0$ (with multiplicity $10$) and 
$12$ (with multiplicity $2$). 
\qed

\medskip

\begin{Le} \la{lepropB}
For all  $n \geq 5$ and $K \subset V$ with $0 \in K$ and $|K| = 5$:
\begin{equation}\la{propB}
B^K(n) := \sum_{J \subset K: |J|=4} \e_J A^{J}(n) \geq 0, 
\end{equation}
where $\e_J$ is the sign of $-c_J$, i.e.\ $\e_J = 1$ if $0 \in J$ 
and $\e_J=-1$ if $0 \notin J$. 
\end{Le}

\proof 
The structure of the proof is very similar to the one of the preceding 
lemma: We consider the block structure of $B^K(n)$ 
corresponding to the $n!/5!$ cosets of $\cX_K^+$ in $\cX_V^+$. 
The diagonal block corresponding to the coset $\cX_K^+$ can be 
identified with $B := B^{\{0,1,2,3,4\}}(5)$
(if $K$ is identified with $\{0,1,2,3,4\}$), and as in the proof 
of Lemma~\ref{lepropA} we see that all 
diagonal blocks are equal and thus copies of $B$ and 
all non-diagonal blocks are $0$. Because of this block decomposition of $B^K(n)$ we only have to show $B \geq 0$. 

$B = B^{\{0,1,2,3,4\}}(5)$ is a symmetric $60 \times 60$ matrix 
with small integer entries that can be computed from \eqref{defmat}. 
Using a computer algebra program one can check that $B$ has the eigenvalues 0 (with multiplicity 45) and 24 
(with multiplicity 15), which implies the assertion for $n=5$. 
However, the following argument allows us to obtain 
the spectrum of $B$ without using a computer. Using
the shorthand notation introduced before Lemma~\ref{lepropA},  
we observe that 
\begin{equation}\la{B5}
B A^{(0)} = 
(A^{(1)} + A^{(2)} + A^{(3)} + A^{(4)} -A^{(0)}) A^{(0)}  = 0.  
\end{equation}
Before proving \eqref{B5} we will use it to compute the spectrum of $B$. 
Let 
$$
B^+ := A^{(1)} + A^{(2)} + A^{(3)} + A^{(4)}  + A^{(0)} \quad 
\text{ so that } \quad  B = B^+ - 2A^{(0)}. 
$$
As an immediate consequence of \eqref{B5}, $(B^+ - 2A^{(0)}) A^{(0)} = 0$ 
and  
$$
A^{(0)}(B^+ - 2A^{(0)})= [(B^+ - 2A^{(0)}) A^{(0)}]^t = 0, 
$$
i.e.\ 
$$
B^+ A^{(0)}= 2(A^{(0)})^2 \quad  \text{ and } \quad 
A^{(0)} B^+= 2(A^{(0)})^2 .
$$
By symmetry we get the same relations for $A^{(i)}$ instead of $A^{(0)}$, 
and by the proof of Lemma~\ref{lepropA} $A^{(i)}$ is a symmetric matrix 
with eigenvalues $0$ and $12$ only, so $(A^{(i)})^2 = 12 A^{(i)}$. 
Using all of these relations we get 
$$
(B^+)^2 = \sum_{i=0}^4 B^+ A^{(i)}= \sum_{i=0}^4 2(A^{(i)})^2
= 24  \sum_{i=0}^4 A^{(i)} = 24 B^{+}
$$
and 
\[
\begin{split}
B^2 &= (B^+ - 2A^{(0)})^2 = (B^+)^2 - 2 B^+ A^{(0)} - 2 A^{(0)} B^+ 
+ 4 (A^{(0)})^2\\ 
&= 24 B^+ - 8 (A^{(0)})^2 + 4 (A^{(0)})^2 
= 24 B^+ - 48 A^{(0)} = 24 B,
\end{split}
\]
i.e.\ $\frac 1 {24} B$ is a projection and thus has eigenvalues $0,1$ only. 
So $B$ has eigenvalues $0,24$ only. Since the trace of $B$ is 
$60 \times (2+2+2+2-2) = 360$, the multiplicity of the eigenvalue 24 
has to be $\frac{360}{24} = 15$, and the multiplicity of the eigenvalue 0 
has to be $45$. 

We will now prove \eqref{B5}, i.e.\ $B=0$ on the image of $A^{(0)}$. 
By the proof of Lemma~\ref{lepropA}, we know the block structure of $A^{(0)}$ 
corresponding to the cosets of $\cX_{\{1,2,3,4\}}^+$ in 
$\cX_{\{0,1,2,3,4\}}^+$: The non-diagonal 
blocks are $0$ and the diagonal blocks are copies of $A$, and 
by \eqref{A4}, the image of $A$ is 
\[
\{a 1_{\eta H} + a' 1_{\eta' H} + a'' 1_{\eta'' H}\,
:\; a,a',a'' \in \bbR\,,\; \text{ such that }\; a+a'+a'' = 0\},
\]
where $H = \{id, (01)(23),(02)(13),(03)(12)\}$ and 
$\eta H$, $\eta' H$, $\eta'' H$ are the three distinct cosets of $H$
in $\cX_{\{0,1,2,3\}}^+$. As usual $1_U$ denotes the indicator function of a given set $U$; 
e.g. $1_{\eta H}$ is the function on $\cX_{\{0,1,2,3\}}^+$ that takes the value $1$ on 
$\eta H$ and the value $0$ otherwise. 
So in particular
\[
\text{Im}(A) \subset \text{Span}(1_{\eta H}: \eta \in \cX_{\{0,1,2,3\}}^+),
\]
and by the block structure of $A^{(0)}$ this implies 
\[
\text{Im}(A^{(0)}) \subset \text{Span}(1_{\eta H^{(0)}}: \eta \in \cX_{\{0,1,2,3,4\}}^+),
\]
where $H^{(0)} = \{id, (12)(34),(13)(24),(14)(23)\}$, and thus 
it suffices to show that $Bv=0$ for every vector $v$ of the form 
$v=1_{\eta H^{(0)}}$, i.e.\ 
\begin{equation} \la{sumB}
\sum_{\s \in H^{(0)}} B_{\eta,\eta'\s} = 0 \qquad \text{ for all } \;
\eta,\eta' \in \cX_{\{0,1,2,3,4\}}^+.
\end{equation}
Since $B_{\eta,\eta'\s}$ only depends on $\eta^{-1}\eta'\s$, 
for the proof of \eqref{sumB} we may assume without loss of generality 
that $\eta = id$. The following observations help to reduce the number 
of choices of $\eta'$ that have to be considered: 
Since $\eta'$ has to be an even permutation of $\{0,1,2,3,4\}$, $\eta'$ has to be $id$, a 3-cycle, a 5-cycle or a product of two disjoint 2-cycles. 
Every 5-cycle necessarily has an entry 0, and the 3-cycle and the 
product of the 2-cycles may or may not contain an entry $0$. 
This gives 6 cases altogether. Since  
$B = A^{(1)} + \ldots + A^{(4)} - A^{(0)}$ and $H^{(0)}$ are 
invariant under permutations of 1,2,3,4, and in each of the above cases 
the permutations differ only by permuting the roles of 1,2,3,4, 
it is sufficient to consider one permutation from each case, say 
$\eta' \in \{id,  (123), (012), (03142), (12)(34),(02)(34)\}$.
Since $(12)(34) \in  id H^{(0)}$ and $(03142),(02)(34) \in 
(012) H^{(0)}$ (see below),
we are done once we  check $ 
\sum_{\s \in H^{(0)}} B_{id,\eta'\s} =0$ for 
$\eta' \in \{id,(123),(012)\}$.
In each of the three cases we compute $\eta' H^{(0)}$ and 
check $f(\eta') := \sum_{\s \in \eta'H^{(0)}} B_{id,\s} =0$ by considering 
the contributions for a fixed $\s$ from $A^{(i)}_{id,\s}$ 
for $1 \leq i \leq 4$ and from $-A^{(0)}_{id,\s}$.\\
(a) If $\eta' = id$, we have $\eta' H^{(0)} = H^{(0)} = \{id, (12)(34),(13)(24),(14)(23)\}$.  $id$ gives a contribution 
of $2+2+2+2-2 = 6$, the others permutations give a contribution of $-2$ each, 
so $f(\eta') = 6-2-2-2=0$. \\
(b) If $\eta' = (123)$, we have $\eta' H^{(0)} = \{(123),(243),(142),(134)\}$, 
and each of these gives a contribution of $1-1=0$. (E.g.\ $(123)$ is a 
3-cycle with entries from $J = \{0,1,2,3\}$ or from $J=\{1,2,3,4\}$.)
Thus $f(\eta') = 0$.\\
(c) If $\eta' = (012)$, we have $\eta' H^{(0)} = \{(012),(02)(34), 
(03142),(04132)\}$. $(012)$ gives the contribution $-1-1$, 
$(02)(34)$ gives 2 and the 5-cycles do not contribute. 
so $f(\eta') = -2 + 2 = 0$. 
\qed

\medskip

\subsection{Proof of Theorem \ref{leoctopus}}
In Subsection~\ref{secred} we have seen that Theorem \ref{leoctopus}
follows once we have shown that the correction matrix is positive 
semi-definite. This can now be obtained from the results of Subsection~\ref{seccor} 
concerning the structure of the correction matrix: 

\begin{Le} \la{corpos}
For every $n \geq 2$ we have $C'(n) \geq 0$.  
\end{Le} 

\proof 
We already have seen that $C'(n) = 0$ for $n=2,3$ and 
$C'(4) = -c_0c_1c_2c_3A \geq 0$. For $n \geq 5$ 
we use the variable $J$ for a subset $J \subset V$ with $|J| = 4$ and 
$K$ for a subset $K \subset V$ with $0 \in K$ and  $|K| = 5$. 
%We are done once we have shown 
The lemma follows from the two inequalities in 
\begin{align*}
0 &\leq \sum_K \frac{|c_K|}{|c_0|} \sum_{J \subset K} \e_J A^{J}(n) 
= \sum_J \Big( \sum_{K \supset J} \frac{\e_J|c_K|}{|c_0|} \Big)  A^{J}(n) 
\\ & \leq \sum_J (-c_J) A^{J}(n) = C'(n).
\end{align*}
The first inequality is an immediate consequence of \eqref{propB}, 
and the second follows from \eqref{propA} once we have checked that
$$
\sum_{K \supset J} \frac{\e_J|c_K|}{|c_0|} \leq -c_J\,,\quad \text{ for all } \;J\,.
$$
If $0 \notin J$, the only set $K \supset J$ containing $0$ is 
$K=J \cup \{0\}$ and we get 
\[
\sum_{K \supset J} \frac{\e_J|c_K|}{|c_0|} 
= \frac{\e_J|c_{J \cup \{0\}}|}{|c_0|} 
= \e_J |c_J| = -c_J. 
\]
If $0 \in J$, the sets $K \supset J$ containing $0$ are
of the form $K = J \cup \{i\}$ with $i \notin J$ and we get 
\[
\sum_{K \supset J} \frac{\e_J|c_K|}{|c_0|} 
= \sum_{i \notin J} \frac{\e_J|c_{J \cup \{i\}}|}{|c_0|} 
= - c_J \frac{\sum_{i \notin J}|c_i|}{|c_0|} \leq -c_J 
\]
since $-c_J \geq 0$ in this case, and 
$$
\sum_{i \notin J} |c_i| 
\leq \sum_{i > 0} |c_i| = \sum_{i > 0} c_i = -c_0 = |c_0|.
$$
\qed

\section{Related Markov chains on weighted graphs}\la{sec_rem}

Here we discuss several stochastic processes that can be associated in a natural way to weighted graphs. 
Each of them is an irreducible, symmetric Markov chain as in Section~\ref{secmc}; in particular 
each one is reversible with respect to the uniform distribution on the corresponding state 
space and has a strictly positive sepectral gap. Furthermore, all of them are sub-processes of the interchange process in the sense of \eqref{mcprojection}, 
which allows us to obtain estimates on their spectral gaps as simple corollaries of Theorem \ref{main}. 
In all examples let $G = (V,E)$ be the complete graph on $n$ vertices -- w.l.o.g. we assume that 
$V = \{1,\ldots,n\}$ -- and $c_{xy} \geq 0$ be given edge weights such that the corresponding skeleton graph 
is connected; see Section~\ref{secprocesses}.

\subsection{Exclusion processes} 

\subsubsection{Symmetric exclusion process}
\begin{figure}[htb]
\psfrag{1}{$1$}
\psfrag{2}{$2$}
\psfrag{3}{$3$}
\psfrag{4}{$4$}
\psfrag{5}{$5$}
\includegraphics[scale = .55]{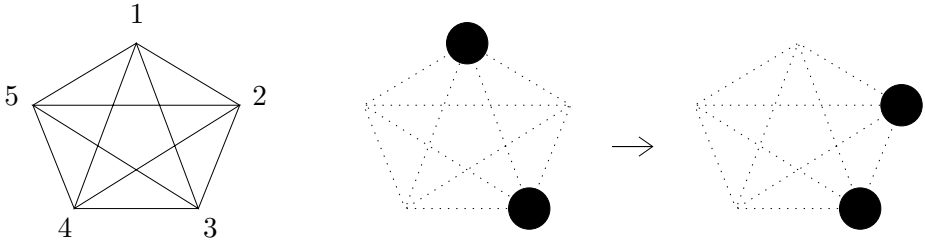}
\caption{2-particle exclusion process on the graph $V=\{1,2,3,4,5\}$. 
The picture shows the underlying graph and a transition
from $\z = \{1,3\}$ to $\z^{1,2} = \{2,3\}$.}  
\label{figexclusion}
\end{figure}

In the $k$-particle exclusion 
process 
a state is an assignment of $k$ indistinguishable particles to 
$k$ of the $n$ vertices of $G$. Here 
$k \in \{1,\ldots,n-1\}$ is a fixed number, which is often omitted in our notation.  
The transition from a state $\z$ to a state 
$\z^{xy}$ (occurring with rate $c_{xy}$) is possible only if in $\z$
one of the positions $x,y$ is occupied and the other is empty. In this 
transition, the particle at the occupied site jumps to the 
empty site; see Figure~\ref{figexclusion}. 
We note that the 1-particle exclusion process is the same as the random walk. 
Formally, the $k$-particle exclusion process is defined to be the Markov chain 
with state space  
$\SEX = \{\z \subset V: |\z| = k\}$ and generator 
$$
\LEX f (\z) = \sum_{xy \in E} c_{xy} (f(\z^{xy}) - f(\z)), \quad 
\text{ where } f: \SEX \to \bbR, \,\z \in \SEX.  
$$
Here $\z^{xy} = \z$ if $xy \subset \z$ or $xy \subset \z^c$ and 
$$
\z^{xy} = \left\{
\begin{aligned}
&(\z \setminus \{y\}) \cup \{x\} 
&& \text{ if $y \in \z$ and $x \notin \z$, }\\
&(\z \setminus \{x\}) \cup \{y\} 
&& \text{ if $x \in \z$ and $y \notin \z$.}
\end{aligned}
\right.
$$
By Section~\ref{secmc}, $-\LEX$ has $|\SEX| = \binom{n}{k}$ 
nonnegative eigenvalues and a positive spectral gap $\lEX_1 >0$.  
The $k$-particle exclusion process can be obtained as a sub-process of the interchange 
process by declaring the sites occupied by particles $1$ through $k$ 
to be occupied and the other vertices to be empty; more precisely
$\pi: \SIC \to \SEX$, $\pi(\eta) = \{\xi_1(\eta),\ldots,\xi_k(\eta)\}$
is a contraction in the sense of \eqref{mcprojection}, which gives 
$\Spec(-\LEX) \subset \Spec(-\LIC)$. In order to compare the 
exclusion process to the random walk,
let $f: V \to \bbR$ be an eigenfunction of 
$-\LRW$ with eigenvalue $\l$ and define $g: \SEX \to \bbR$ by 
$g(\z) = \sum_{x \in \z} f(x)$. Note that if $g$ is constant then 
$f$ must be constant. 
Therefore, $g \not \equiv 0$ (since otherwise $f$ is constant and thus 
$f \equiv 0$), and  $\sum_{x,y \in \z, x\neq y} c_{xy} (f(y) - f(x)) = 0$ implies
\[
\begin{split}
&(-\LEX g)(\z) 
= - \sum_{x \in \z, y \notin \z} c_{xy}(g(\z^{xy})-g(\z))\\
&= - \sum_{x \in \z, y \notin \z} c_{xy}(f(y)-f(x))
= - \sum_{x \in \z, y \neq x} c_{xy}(f(y)-f(x))\\
&= \sum_{x \in \z} (-\LRW f) (x) 
= \l \sum_{x \in \z} f (x) 
= \l g(\z)\,,
\end{split}
\]
i.e.\ $g$ is an eigenfunction of $-\LEX$ with eigenvalue $\l$. This gives 
$$
\Spec(-\LRW) \subset \Spec(-\LEX) \subset \Spec(-\LIC) \,, 
$$
As a corollary of Theorem \ref{main} one has that, for arbitrary 
number of particles 
$k=1,\dots,n-1$, for every graph $G$:
\be\la{cor1}
\lEX_1(G) = \lRW_1(G)
\end{equation}

\subsubsection{Colored exclusion process.} 
\begin{figure}[htb]
\psfrag{1}{$1$}
\psfrag{2}{$2$}
\psfrag{3}{$3$}
\psfrag{4}{$4$}
\psfrag{5}{$5$}
\includegraphics[scale = .55]{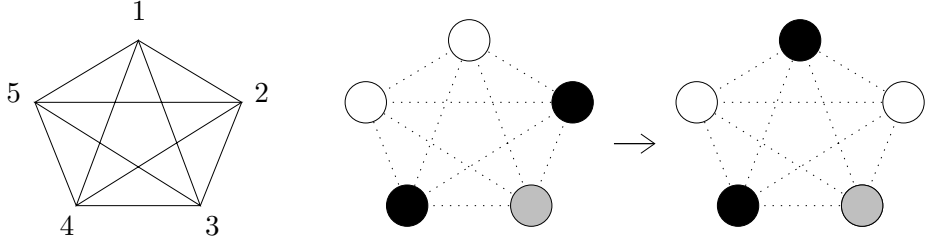}
\caption{Colored exclusion process on $V=\{1,2,3,4,5\}$
with 3 types of particles $(n_1 = 2,n_2 =1, n_3=2)$.  
The picture shows the underlying graph and a transition
from $\a = (\{2,4\},\{3\},\{1,5\})$ to $\a^{1,2} = (\{1,4\},\{3\},\{2,5\})$.}  
\label{figmulti}
\end{figure}

In the colored exclusion process there are $r \geq 2$ types of particles 
($n_i \geq 1$ of type $i$ such that $n_1+ \cdots + n_r = n$), 
where particles of the same type (or color) are indistinguishable. 
A state is an assignment of these particles
to the vertices of $G$ so that every vertex is occupied by exactly one 
particle, and in the transition from a state $\a$ to 
a state $\a^{xy}$ particles at sites $x$ and $y$ interchange 
their positions; see Figure~\ref{figmulti}.
Formally, the colored exclusion process is the Markov chain 
on the state space $\SCEP$, which is the set of partitions 
$\a=(\a_1,\ldots,\a_r)$ of $V$ such that $|\a_i| = n_i$, and 
the generator is defined by 
$$
\LCEP f (\a) = \sum_{xy \in E} c_{xy} (f(\a^{xy}) - f(\a)), \quad 
\;\; f: \SCEP \to \bbR, \a \in \SCEP.  
$$
Here $\a^{xy} = \a$ if $x,y \in \a_i$ for some $i$, and if 
$x \in \a_i$ and $y \in \a_j$ for $i \neq j$ we have 
$\a^{xy} = (\a_1^{xy}, \ldots, \a_r^{xy})$, where 
$\a_i^{xy} = (\a_i \setminus \{x\}) \cup \{y\}$, 
$\a_j^{xy} = (\a_j \setminus \{y\}) \cup \{x\}$, 
and $\a_k^{xy} = \a_k$ for all $k \neq i,j$. 
By Section~\ref{secmc}, $-\LCEP$ has $|\SCEP| = \binom{n}{n_1,\ldots,n_r}$ 
nonnegative eigenvalues and a positive spectral gap $\lCEP_1 >0$.  
The $n_1$-particle exclusion process is a sub-process of the 
colored exclusion process (by declaring all sites occupied by 
particles of type $2,\ldots,r$ to be empty), which in turn 
is a sub-process of the interchange process (by declaring particles 
$1,\ldots,n_1$ to be of type 1, \ldots, particles $n-n_r+1, \ldots, n$ to 
be of type $r$). The definitions of the corresponding contractions are 
obvious. This gives, for the given choice of parameters $n_1,\dots,n_r$, 
$$
\Spec(-\LEX) \subset \Spec(-\LCEP) \subset \Spec(-\LIC)\,.
$$
From Theorem \ref{main} and (\ref{cor1}) 
it follows that for any choice of the parameters $r,n_1,\dots,n_r$,
for any graph $G$: 
\be\la{cor2}
\lCEP_1(G)=\lRW_1(G)\,.
\end{equation}

\subsection{Cycles and matchings.} 
We turn to examples 
of processes with a gap that is in general  
strictly larger than that of the random walk. 
We note that the processes defined here are examples 
from a general class of processes, obtained as 
the evolution of certain subgraphs
of the complete graph when
the labels undergo the interchange process on $G$. 
For all these processes one has analogous
estimates for the spectral gap.

\subsubsection{Cycle process}
\begin{figure}[htb]
\psfrag{1}{$1$}
\psfrag{2}{$2$}
\psfrag{3}{$3$}
\psfrag{4}{$4$}
\psfrag{5}{$5$}
\includegraphics[scale = .55]{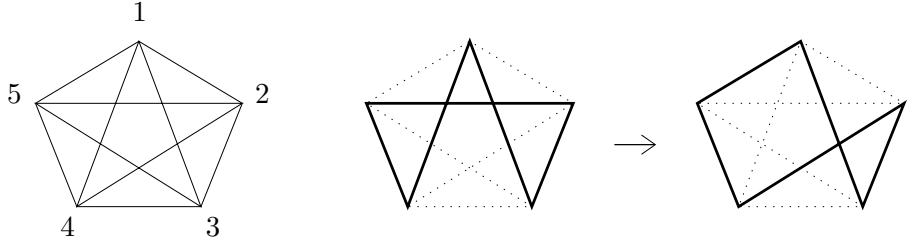}
\caption{Cycle process on $V=\{1,2,3,4,5\}$. 
The picture shows the underlying graph and a transition 
from $\g = \{\{1,3\},\{3,2\},\{2,5\},\{5,4\},\{4,1\}\}$ to 
$\g^{1,2} = \{\{1,3\},\{3,2\},\{2,4\},\{4,5\},\{5,1\}\}$.}  
\label{figcycle}
\end{figure}

The states of 
the cycle process are $n$-cycles, where $n = |V|$. In order to avoid 
a trivial situation we assume $n  \geq 4$. One could think of a rubber band 
that at certain points is pinned to the vertices of $G$. 
The transition from $\g$ to $\g^{xy}$ (occurring with rate $c_{xy}$) 
can be obtained by taking the point of the rubber band pinned to $x$  
from $x$ to $y$ and the point pinned to $y$ from $y$ to $x$; 
see Figure~\ref{figcycle}.
Formally, an $n$-cycle in $G$ is 
a set of edges $\g \subset E$ forming a subgraph of $G$ 
isomorphic to $\{\{1,2\},\ldots ,\{n-1,n\},\{n,1\}\}$. 
The cycle process is the Markov chain with state space $\SCY$, 
the set of all $n$-cycles of $G$, and generator 
$$
\LCY f (\g) = \sum_{xy \in E} c_{xy} (f(\g^{xy}) - f(\g)), \quad 
\text{ where } f: \SCY \to \bbR, \,\g \in \SCY.  
$$
Here $\g^{xy} = \{ b^{xy}: b \in \g\}$, where for an edge $b \in E$ 
we define $b^{xy} = b$ if $x,y \notin b$ or $b = xy$, 
and $b^{xy} = xz$ if $b=yz$ ($z \neq x,y$). 
%By Section~\ref{secmc}
Thus, $-\LCY$ has $|\SCY| = \frac{(n-1)!}{2}$ 
nonnegative eigenvalues and a positive spectral gap $\lCY_1 >0$.  
The cycle process can be obtained from the interchange process by pinning 
a cycle on the particles labeled $1,\ldots,n,1$ in that order, i.e.\ 
$\pi: \SIC \to \SCY$, $\pi(\eta) = \{\xi_1(\eta)\xi_2(\eta), \ldots, 
\xi_{n-1}(\eta)\xi_n(\eta),\xi_n(\eta)\xi_1(\eta)\}$ is a contraction. 
This gives 
$
\Spec(-\LCY) \subset \Spec(-\LIC)$, %\quad \text{ and thus }
and thus, by Theorem \ref{main}, on any graph $G$:
\begin{equation}\la{cor3}
\lCY_1(G) \geq \lRW_1(G)\,.
\end{equation}
It is possible to see that in general the inequality is strict, e.g.\ 
for $n=4$ and $c_{ij} = 1$ for all $ij \in E$ one can compute
$\lRW_1=4$ and $\lCY_1=6$. 

%So Conjecture~\ref{aldous} would imply $\lCY_1 \geq \lRW_1$. 

\subsubsection{Matching process}
\begin{figure}[htb]
\psfrag{1}{$1$}
\psfrag{2}{$2$}
\psfrag{3}{$3$}
\psfrag{4}{$4$}
\psfrag{5}{$5$}
\psfrag{6}{$6$}
\includegraphics[scale = .55]{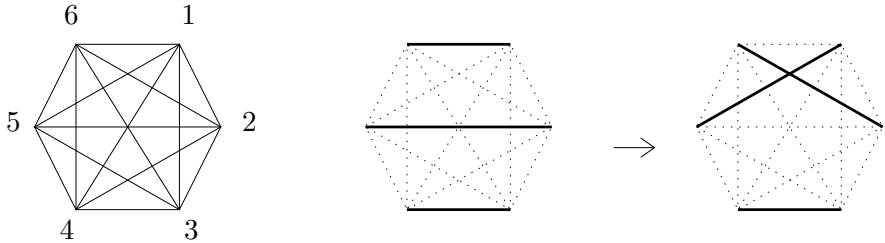}
\caption{Matching process on $V=\{1,2,3,4,5,6\}$. 
The picture shows the underlying graph and a transition 
from $\zeta = \{\{1,6\},\{2,5\},\{3,4\}\}$ to 
$\zeta^{1,2} = \{\{2,6\},\{1,5\},\{3,4\}\}$.}  
\label{figmatch}
\end{figure}

Set 
$n=2k$, $k \geq 2$. A perfect matching of $n$ elements is given by $k$ disjoint 
edges of the complete graph on $V=\{1,\dots,n\}$. A perfect matching
configuration is thus denoted $\zeta=\{b_1,\dots,b_k\}$, where the edges 
$b_i\in E$ are such that $b_i\cap b_j =\emptyset$. 
The state space of the matching process, denoted $\SMP$, is the set of all possible perfect
matchings $\zeta$. The transition $\zeta\to\zeta^{xy}$, occurring with
rate $c_{xy}$, is described as follows: if $xy\in\zeta$ then nothing
happens, and $\zeta^{xy}=\zeta$; if $xy\notin \zeta$ then there are
uniquely determined 
points $u,v\in V$ 
such that $xu,yv\in \zeta$, and $\zeta^{xy}$ coincides with
$\zeta$ except that $xu,yv$ are removed while $yu,xv$ are added; 
see Figure~\ref{figmatch}.
The generator is then given by 
$$
\LMP f (\zeta) = \sum_{xy \in E} c_{xy} (f(\zeta^{xy}) - f(\zeta)), \quad \text{ where }
f: \SMP \to \bbR, \,\zeta \in \SMP\,.  
$$
Thus, $-\LMP$ has $|\SMP| = \frac{(2k)!}{2^k k!}$ 
nonnegative eigenvalues and a positive spectral gap $\lMP >0$.  
The matching process is a sub-process of the interchange process in
the sense of (\ref{mcprojection}). Indeed, given a permutation
$\eta\in\cX_n$, one can obtain a perfect matching by setting 
$$
\zeta = \pi(\eta) = \{\xi_1(\eta)\xi_{k+1}(\eta),\dots,\xi_{k}(\eta)\xi_{2k}(\eta)\}\,,
$$
and the map $\pi: \SIC \to \SMP$ defines the desired contraction. 
This shows that $
\Spec(-\LMP) \subset \Spec(-\LIC)$, %\quad \text{ and thus }
which implies, by Theorem \ref{main}, that for all graphs $G$:
\begin{equation}\la{cor4}
\lMP_1(G) \geq \lRW_1(G)\,.
\end{equation}
It is known that for the unweighted complete graph $c_{xy}\equiv 1$ the inequality (\ref{cor4}) is strict. 
We refer to \cite{DH} for a complete description of $\Spec(-\LMP)$ in this special case, and note that 
for $k=2,n=4$ the two eigenvalues are 6 and 4 respectively.

\bigskip
\bigskip
 
\noindent
{\bf Acknowledgments}.  P.C.\ thanks Filippo Cesi for helpful discussions. 
Partial support from  
NSF Grant DMS-0301795 is acknowledged. 
P.C.\ was also partially supported by the 
Advanced Research Grant ``PTRELSS'' ADG-228032 of the European Research Council.

\end{document}